\documentclass[a4paper,10pt]{article}

\usepackage{amsfonts,amssymb,amsmath,xypic}

\numberwithin{equation}{section}

\usepackage{enumitem}

\usepackage{tikz-cd}

\newtheorem{theo}[equation]{Theorem}
\newtheorem{coro}[equation]{Corollary}
\newtheorem{lemm}[equation]{Lemma}
\newtheorem{prop}[equation]{Proposition}
\newtheorem{defi}[equation]{Definition}
\newtheorem{rema}[equation]{Remark}
\newtheorem{exam}[equation]{Example}

\newenvironment{proof}{\noindent \textbf{{Proof.}} \sf}

\def\qed{\hfill $\diamond$ \bigskip}

\def\C{{\mathcal C}}

\def\s{{\mathcal S}}

\def\lim{\mathop{\rm lim}\nolimits}

\def\HH{\mathsf H}
\def\HHH{\mathsf{HH}}

\def\Ext{\mathsf{Ext}}
\def\Hom{\mathsf{Hom}}
\def\Tor{\mathsf{Tor}}

\def\Ker{\mathsf{Ker}}

\def\Coker{\mathsf{Coker}}

\def\dim{\mathsf{dim}}

\def\tauh{\HHH^1_\tau(\Lambda)}
\def\sfD{\mathsf{D}}
\def\sfTr{\mathsf{Tr}}

\begin{document}
\sf

\title{On the first $\tau$-Hochschild cohomology of an algebra}
\author{Claude Cibils,  Marcelo Lanzilotta, Eduardo N. Marcos,\\and Andrea Solotar
\thanks{\tiny This work has been supported by the projects PIP-CONICET 11220200101855CO, USP-COFECUB.
The third mentioned author was supported by the thematic project of FAPESP  2018/23690-6,  research grants from CNPq 308706/2021-8  and  310651/2022-0. The fourth mentioned author is a research member of CONICET (Argentina), Senior Associate at ICTP and visiting Professor at Guangdong Technion-Israel Institute of Technology.}}

\date{}
\maketitle
\begin{abstract}

  In this paper we introduce, according to one of the main ideas of $\tau$-tilting theory,  the $\tau$-Hochschild cohomology in degree one of a finite dimensional  $k$-algebra $\Lambda$, where $k$ is a field. We define the excess of $\Lambda$  as the difference between the  dimensions of the $\tau$-Hochschild  cohomology in degree one and  the dimension of the usual Hochschild cohomology in degree one.

 One of the main results is that for a zero excess bound quiver algebra $\Lambda=kQ/I$, the Hochschild cohomology in degree two $\HHH^2(\Lambda) $ is isomorphic to the space of morphisms $\Hom_{kQ-kQ}(I/I^2, \Lambda).$ This is useful to determine when $\HHH^2(\Lambda)=0$ for these algebras.

We compute the excess for hereditary, radical square zero and  monomial triangular algebras. For a bound quiver algebra $\Lambda$, a formula for the excess of $\Lambda$ is obtained.  We also give  a criterion for $\Lambda$ to be $\tau$-rigid.

\end{abstract}




\section{\sf Introduction}

Let $A$ be a finite dimensional algebra over a field $k$,   that  we will call   an algebra for short. Let $M$ and $N$ be finitely generated left $A$-modules, henceforth called left $A$-modules. Let $\tau$  denote  the Auslander-Reiten translation, see for instance \cite{AUSLANDERREITENSMALO1995} or \cite{IYAMAREITEN2014},  and denote $\sfD(-)=\Hom_k(-,k)$.  We reproduce an extract from B. Marsh's lecture notes in Cologne \cite{MARSH2023}: ``the Auslander-Reiten duality suggests that in contexts where $\Ext^1_A (M,N)$ appears, we might investigate replacing it with $\sfD\Hom_A(N,\tau M)$ and this can be regarded as one of the main ideas of $\tau$-tilting theory." While $\sfD$ is absent in  the original text, $\sfD$ is present  in Auslander-Reiten's duality formula for it to be functorial. Of course  adding $\sfD$ does not change the dimensions. Recall that $M$ is called $\tau$-\emph{rigid} if $\Hom_A(M,\tau M)=0$, see for instance \cite[Subsection 4.1]{IYAMAREITEN2014}.

 On the other hand, let $\Lambda^e=\Lambda\otimes_k \Lambda^{\mathsf{ op}}$ be the enveloping algebra of an algebra $\Lambda$. Let $X$ be a $\Lambda$-bimodule. The Hochschild cohomology of $\Lambda$ with coefficients in $X$ is $\HH^n(\Lambda, X)= \Ext^n_{\Lambda^e}(\Lambda, X)$, see \cite{CARTANEILENBERG, HOCHSCHILD1945,WITHERSPOON2019} and it is denoted  $\HHH^n(\Lambda)$ when $X=\Lambda$.  Moreover, Hochschild homology is $\HH_n(\Lambda,X)= \Tor_n^{\Lambda^e}(\Lambda,X)$.   Since left $\Lambda^e$-modules are the same as $\Lambda$-bimodules, in the sequel we often replace $\Lambda^e$ with $\Lambda - \Lambda$.

According to the main idea of  $\tau$-tilting theory  mentioned above, we will investigate in this paper   the replacement of    $\Ext^1_{\Lambda-\Lambda}(\Lambda,X)$ by  the  $\tau$-\emph{Hochschild cohomology in degree one} $\HH^1_\tau(\Lambda,X)= \sfD\Hom_{\Lambda-\Lambda} (X, \tau \Lambda)$.  Note that  here $\tau$ is the Auslander-Reiten translation of left  $\Lambda^e$-modules, that is of $\Lambda$-bimodules. When  $X=\Lambda$, we denote $\tauh=\HH^1_\tau(\Lambda,\Lambda)$. The excess $e(\Lambda)$ is defined as $\dim_k \tauh - \dim_k \HHH^1(\Lambda)$.

 One of  the main results of this paper is that for a zero excess bound quiver algebra $\Lambda=kQ/I$ we have $\HHH^2(\Lambda) = \Hom_{kQ-kQ}(I/I^2, \Lambda)$ -- see Corollary \ref{H2 is zero}. This result will be useful in a future work to determine when an algebra with zero excess has zero Hochschild cohomology in degree $2$. The algebras $\Lambda$ with $\HHH^2(\Lambda)=0$ are important since they are rigid in the following sense. Suppose that $k$ is algebraically closed and let $V$ be a $k$-vector space of dimension $n$. Let ${\mathcal Alg}_n$ be the affine open subscheme of algebra structures with $1$ of the affine algebraic scheme defined by
   ${\mathcal S}_n(R)=\{\mbox{associative $R$-algebra structures on } R\otimes_k V\},$  where $R$ is a commutative $k$-algebra.  Corollary 2.5 of \cite{GABRIEL1975} states that $\HHH^2(\Lambda)=0$ if and only if the orbit of $\Lambda\in {\mathcal Alg}_n$ under the general linear group ${\mathcal GL}(V)$ is an open subscheme of ${\mathcal Alg}_n$ -- that is by definition, $\Lambda$ is rigid. Moreover, P. Gabriel in \cite[p. 140]{GABRIEL1975} mentions that it should be one of the main tasks of associative algebra to determine for every $n$ the number of irreducible components of ${\mathcal Alg}_n$. The determination of open orbits makes it possible to obtain lower bounds for the number of irreducible components of ${\mathcal Alg}_n$, as G. Mazzola did in \cite[p. 100]{MAZZOLA 1979}.

The paper is organised as follows. In Section \ref{tau-tilting H} we give a more detailed definition, as well as properties of the $\tau$-Hochschild cohomology and of the excess.   Let $\mathsf{Tr}$ be the transpose of a bimodule, see for instance \cite{AUSLANDERREITENSMALO1995}, and recall that $X_\Lambda$ denotes the coinvariants of a $\Lambda$-bimodule $X$, see Remark \ref{invariants and coinvariants}.   We prove that $\tauh= \mathsf{(Tr \Lambda})_\Lambda$ and we give a formula for the dimension of the vector space $\tauh$.  Note that in this article we use the symbol = for the existence of a canonical isomorphism.

In Section \ref{cases}, for an hereditary algebra $\Lambda$ we  prove that the dimensions of $\tauh$ and $\HHH^1(\Lambda)$ are equal.
We say that an algebra $\Lambda$ has the  $\HH^2$ \emph{cancellation properties}  if $\HHH^2(\Lambda)=0=\HH^2(\Lambda, r^i)$ for all $i>0$, where $r$ is the Jacobson radical of $\Lambda$. For instance hereditary algebras have the $\HH^2$ cancellation properties. We obtain that $e(\Lambda)=0$ whenever $\Lambda$ has the $\HH^2$ cancellation properties, based on a formula for the dimension of $\HHH^1(\Lambda)$ in \cite{CIBILS1988}.

In Section \ref{cases} we  also consider radical square zero algebras and monomial algebras whose quiver has no oriented cycles. For those algebras $\Lambda$ we prove that $\HHH^1(\Lambda)=0$ if and only if $\tauh =0$, and this occurs precisely when $Q$ is a tree. This extends a result of \cite{BARDZELL MARCOS 1998}. We provide examples where   the excess is not zero.

In Theorem \ref{exceso} we give a formula for the excess of a bound quiver algebra. Finally we provide a criterium for the algebra $\Lambda$ to be $\tau$-rigid in terms of the dimension of its Hochschild cohomology in degree $2$.

In a future article, we will consider higher $\tau$-Hochschild cohomology in relation to Happel's question, see \cite{HAPPEL1988}.

\section{\sf $\tau$-Hochschild cohomology in degree one}\label{tau-tilting H}

We begin this section by briefly recalling the definition of the Auslander-Reiten translation and the duality formula which is useful for our aims, for more details see for instance \cite{AUSLANDERREITENSMALO1995} or \cite{ANGELERI2006}. Let $A$ be an algebra and $M$  a left $A$-module.

First, the transpose $\sfTr M$ is defined as follows. Consider a minimal projective presentation of $M$
$$P_1 \stackrel{d_1}{\longrightarrow} P_0 \stackrel{d_0}{\longrightarrow}M\longrightarrow 0.$$
 Applying  to $d_1$ the functor $\Hom_A(-, A)$ which sends left $A$-modules to right $A$-modules we get
$$\Hom_A (P_0, A) \stackrel{d_1^*}{\longrightarrow} \Hom_A (P_1, A).$$
By definition, $\sfTr M = \mathsf{Coker} d_1^*$.

This gives a bijection between the isomorphism classes  of indecomposable non-projective left $A$-modules and the isomorphism classes  of indecomposable non-projective right $A$-modules.

Next, the  exact  functor $\sfD =\Hom_k(-,k)$ sends right $A$-modules to left $A$-modules. We obtain an exact sequence of left $A$-modules
$$0 \longrightarrow \sfD\sfTr M \longrightarrow \sfD\Hom_A (P_1, A) \stackrel{\sfD d_1^*}{\longrightarrow} \sfD\Hom_A (P_0, A).$$

 Finally  by definition $\tau M= \sfD\sfTr M$.

This gives a bijection between the isomorphism classes of indecomposable non-projective left $A$-modules and the isomorphism classes of indecomposable non-injective left $A$-modules.

Let $M$ and $N$ be left $A$-modules. Let  $\mathcal{I}  \Hom_A(M,N)$  be the  $k$-subspace  of $\Hom_A(M,N)$ of morphisms which factor through an injective left $A$-module. The quotient is denoted $\overline{\Hom}_A(M,N)$. The Auslander-Reiten duality formula in \cite{AUSLANDERREITEN1975} is
$$\Ext_A^1(M,N) =    \sfD\overline{\Hom}_A(N,\tau M).$$

As mentioned in the Introduction, one of the main ideas of $\tau$-tilting theory is to replace $\Ext^1_A (M,N)$ with $\sfD\Hom_A(N,\tau M)$, which in a sense amounts to recover the missing morphisms which factor  through  injectives.

Let $\Lambda$ be an algebra. To define  the  $\tau$-Hochschild cohomology in degree one, recall that $\HH^1(\Lambda,X) = \Ext^1_{\Lambda-\Lambda}(\Lambda,X)$. Note that this concerns bimodules, hence in the following $\tau$ is the Auslander-Reiten translation for bimodules  or  equivalently of left $\Lambda^e$-modules.

\begin{defi}
Let $\Lambda$ be an algebra, and let $X$ be a $\Lambda$-bimodule. The first \emph{$\tau$-Hochschild cohomology}  of $\Lambda$  with coefficients in $X$ is
$$\HH^1_\tau(\Lambda,X)= \sfD\Hom_{\Lambda -\Lambda}(X, \tau \Lambda).$$
\end{defi}

In this paper we will focus in the case $X=\Lambda$: $$\tauh=\sfD\Hom_{\Lambda -\Lambda}(\Lambda, \tau \Lambda).$$

\begin{defi}\label{excess}
The \emph{excess} of an algebra $\Lambda$ is
$$e(\Lambda) = \dim_k \tauh - \dim_k \HHH^1(\Lambda).$$
\end{defi}

\begin{lemm}\label{excess is dim of factoring injectives}
 The excess is a non negative integer, equal to $\dim_k \mathcal{I}\Hom_{\Lambda-\Lambda} (\Lambda, \tau\Lambda).$
\end{lemm}

\begin{proof}
By definition $$\tauh=\sfD\Hom_{\Lambda -\Lambda}(\Lambda, \tau \Lambda),$$ while $$\HHH^1(\Lambda)= \Ext^1_{\Lambda -\Lambda}(\Lambda,\Lambda).$$ The Auslander-Reiten duality formula is
$$\Ext_{\Lambda -\Lambda}^1(\Lambda,\Lambda) =    \sfD\left(\frac{\Hom_{\Lambda-\Lambda}(\Lambda ,\tau \Lambda)}{\mathcal{I}\Hom_{\Lambda-\Lambda}(\Lambda ,\tau \Lambda)}\right).$$\qed
\end{proof}

Next we recall some well known facts about invariants and coinvariants.

\begin{rema}\label{invariants and coinvariants}
Let $\Lambda$ be an algebra and let $X$ be a $\Lambda$-bimodule.
\begin{itemize}
  \item
The subspace of invariants  of $X$  is
$$\HH^0(\Lambda,X)= X^\Lambda=\{x\in X \ | \ \forall \Lambda \in \Lambda, \Lambda x= x\Lambda\}= \Hom_{\Lambda-\Lambda}(\Lambda, X)$$
where the last canonical isomorphism sends $\varphi \in \Hom_{\Lambda-\Lambda}(\Lambda, X)$ to $\varphi(1)$.
\item
The vector space of coinvariants   of $X$   is
$$\HH_0(\Lambda, X) = X_\Lambda= X/\langle\lambda x - x\lambda \ |\ \lambda\in\Lambda \mbox{ and } x\in X\rangle = \Lambda\otimes_{\Lambda-\Lambda} X$$
where the last canonical isomorphism sends $\lambda\otimes x\in \Lambda\otimes_{\Lambda-\Lambda} X$ to the class of $\lambda x$.
\item
It is easy to show that $D\left(X^\Lambda\right)=(DX)_\Lambda$. Observe that more generally we have   in all degrees
$$DH^n(\Lambda, X) = \HH_n(\Lambda, DX).$$

\end{itemize}

\end{rema}
\begin{prop}\label{tau tilting cohomology is coinvariants of the transpose}
Let $\Lambda$ be an algebra. We have
$$\tauh =(\sfTr \Lambda)_\Lambda.$$
\end{prop}

\begin{proof}
 Let $Y=\sfD\sfTr\Lambda$. According to Remark \ref{invariants and coinvariants} and using that   $\sfD^2$ is the identity, we have the following chain of equalities and canonical isomorphisms of vector spaces:
$$\tauh=\sfD \Hom_{\Lambda-\Lambda}(\Lambda, Y) =\sfD \left(Y^\Lambda\right)= (\sfD Y)_\Lambda = (\sfD\sfD\sfTr\Lambda)_\Lambda=(\sfTr\Lambda)_\Lambda.$$\qed
\end{proof}

In this paper a \emph{quiver} $Q$ is a finite oriented graph, given by a set of vertices $Q_0$, a set of arrows $Q_1$, and two maps called source and target $s,t : Q_1 \to Q_0$. The \emph{quiver algebra} $kQ$ is a vector space with basis the set $B$ of all oriented paths in $Q$, including those of length $0$,  that is $Q_0$. The  product of two paths is their concatenation if it is possible and $0$ otherwise. The algebra structure of $kQ$ is obtained by extending linearly the  product on paths. Note that $Q_0$ is a set of orthogonal idempotents, their sum gives the unit of $kQ$. The set of paths of strictly positive length $B^{>0}$ is a basis of the ideal $F=\langle Q_1 \rangle$. An ideal $I$ is \emph{admissible} if there exists $n\geq 2$ such that  $F^n \subset I \subset F^2$. The quotient algebra $kQ/I$ is called a \emph{bound quiver algebra.}

An algebra $\Lambda$ is called \emph{sober}  if the endomorphism algebra of each simple left $\Lambda$-module is reduced to $k$, which is always the case if $k$ is algebraically closed. A well known result of P.Gabriel is that any sober algebra is Morita equivalent to a bound quiver algebra $kQ/I$ for a unique quiver $Q$. Note that the admissible ideal $I$ is in general not unique.

\begin{theo}\label{dimension of tau tilting cohomology}
Let $\Lambda=kQ/I$ a bound quiver algebra, and let $Z\Lambda$ be its center. We have
$$\dim_k\tauh= \dim_kZ\Lambda -\sum_{x\in Q_0}\dim_k x\Lambda x + \sum_{a\in Q_1} \dim_k t(a)\Lambda s(a).$$
\end{theo}
\begin{proof} By Proposition \ref{tau tilting cohomology is coinvariants of the transpose}, we have to compute $\dim_k (\sfTr \Lambda)_\Lambda$. To begin with, we will consider $\sfTr \Lambda$. Let $E=kQ_0$, which is a maximal commutative semisimple subalgebra of $kQ$. The projective minimal presentation of $\Lambda$ as $\Lambda$-bimodule is known to have the following form, see \cite[p. 324]{BUTLER KING 1997} and \cite[p. 72]{BARDZELL 1997}
\begin{equation}\label{projective presentation of Lambda}
\Lambda\otimes_E kQ_1 \otimes_E \Lambda \stackrel{f}{\longrightarrow} \Lambda\otimes_E \Lambda \longrightarrow \Lambda \longrightarrow  0
\end{equation}
where $\Lambda\otimes_E \Lambda \longrightarrow \Lambda$ is given by the product of $\Lambda$. For $a\in Q_1$ we have
$$f(t(a)\otimes a\otimes s(a))= a\otimes s(a) - t(a)\otimes a.$$
Consequently, for $\lambda$, $\mu\in \Lambda$ we obtain
$$f(\mu\otimes a\otimes \lambda)= \mu a\otimes s(a)\lambda - \mu t(a)\otimes a\lambda.$$
We write $\otimes$ instead of $\otimes_k$. Also note that the enveloping algebra $\Lambda^e$ viewed as a $\Lambda$-bimodule is isomorphic to $\Lambda\otimes \Lambda$ with action $\lambda(a\otimes b)\mu = \lambda a\otimes b\mu.$

The functor $\Hom_{\Lambda-\Lambda}(- , \Lambda\otimes\Lambda)$ applied to (\ref{projective presentation of Lambda})  provides the  exact sequence defining $\sfTr\Lambda$
$$ \Hom_{\Lambda-\Lambda}(\Lambda\otimes_E \Lambda, \Lambda\otimes\Lambda) \stackrel{f^*}{\longrightarrow} \Hom_{\Lambda-\Lambda}( \Lambda\otimes_E kQ_1 \otimes_E \Lambda , \Lambda\otimes\Lambda) \longrightarrow \sfTr \Lambda\longrightarrow 0.$$
   Next we use   that for an $E$-bimodule $U$ and a $\Lambda$-bimodule $X$   there is a canonical isomorphism
$$  \Hom_{\Lambda-\Lambda}( \Lambda\otimes_E U \otimes_E \Lambda , X)=  \Hom_{E-E}(U,X)$$
and observe that $\Lambda\otimes_E \Lambda=\Lambda\otimes_E E \otimes_E\Lambda$.    We thus obtain the following exact sequence, where we kept the same notation  for the $\Lambda$-bimodule morphism $f^*$
\begin{equation}\label{exact Tr E-E}
\Hom_{E-E}(E,\Lambda\otimes\Lambda)\stackrel{f^*}{\longrightarrow}\Hom_{E-E}(kQ_1  , \Lambda\otimes\Lambda) \longrightarrow \sfTr\Lambda\longrightarrow 0.
\end{equation}
In the following we work out the  exact sequence (\ref{exact Tr E-E}). Let $y,x\in Q_0$ and let ${}_yk_x$ be the simple $E$-bimodule of dimension 1 given by the idempotent $y\otimes x \in E^e$, namely ${}_yk_x = yE\otimes Ex$. Let $U$ be an $E$-bimodule. Clearly we have   a canonical isomorphism
$$ \Hom_{E-E} ({}_yk_x, U) = yUx.$$
Observe that as $E$-bimodules we have $$E= \oplus_{x\in Q_0}\ {}_xk_x \mbox{\  \  and\ \  } kQ_1=\oplus_{a\in Q_1} \ {}_{t(a)}k_{s(a)}.$$
The exact sequence (\ref{exact Tr E-E}) becomes, by still keeping the same notation for $f^*$
\begin{equation}\label{quasi final Tr}
\oplus_{x\in Q_0} \ (x\Lambda\otimes\Lambda x) \stackrel{f^*}{\longrightarrow} \oplus_{a\in Q_1}\  \left(t(a)\Lambda\otimes \Lambda s(a)\right) \longrightarrow \sfTr\Lambda\longrightarrow 0.
\end{equation}
 Let $M$ be a right $\Lambda$ - module and $N$ be a left $\Lambda$-module,  $M\otimes N$ is a $\Lambda$-bimodule  for the \emph{internal} action $\lambda(m\otimes n)\mu =m\mu\otimes \lambda n$. On the other hand  $N\otimes M$ is a $\Lambda$-bimodule for the \emph{external} action $\lambda(n\otimes m)\mu =\lambda n\otimes  m\mu$. Of course, these $\Lambda$-bimodules are isomorphic  through the flip map $\sigma(n\otimes m)= m\otimes n$.

We rewrite \ref{quasi final Tr} using the flips maps $$\sigma_x : x\Lambda \otimes \Lambda x \to \Lambda x \otimes x\Lambda \mbox{ and } \sigma_a : t(a)\Lambda\otimes \Lambda s(a) \to \Lambda s(a)\otimes t(a)\Lambda$$ thus getting an exact sequence for bimodules  with external action. By abuse of notation   we  still write $f^*$ instead of $\left(\oplus_{a\in Q_1}\sigma_a\right) f^* \left(\oplus_{x\in Q_0}\sigma_x^{-1}\right)$.
\begin{equation}\label{final Tr}
\oplus_{x\in Q_0} (\Lambda x\otimes x\Lambda)  \stackrel{f^*}{\longrightarrow} \oplus_{a\in Q_1} (\Lambda s(a)\otimes t(a)\Lambda) \longrightarrow \sfTr\Lambda\longrightarrow 0.
\end{equation}

It is an easy but rather meticulous computation to   track  the morphism of $\Lambda$-bimodules $f^*$ along the   previous  steps.  In the end, we obtain the following formula in the context of (\ref{final Tr}):
\begin{equation}\label{f^*}
 f^*(x\otimes x)= \sum_{ \substack{a\in Q_1\\s(a)=x}} x\otimes a \ \ -\ \sum_{ \substack{b\in Q_1\\t(b)=x}} b\otimes x.
\end{equation}
Recall that our aim is to compute the dimension of the coinvariants of $\sfTr\Lambda$,   that is of $\Lambda\otimes_{\Lambda-\Lambda} \sfTr \Lambda$ by Remark \ref{invariants and coinvariants}.  The functor $\Lambda\otimes_{\Lambda-\Lambda} -$ is right exact and preserves direct sums, so we obtain the exact sequence
\begin{equation}\label{coinvariants Tr}
\oplus_{x\in Q_0} (\Lambda x\otimes x\Lambda)_\Lambda  \stackrel{f^*_\Lambda}{\longrightarrow} \oplus_{a\in Q_1} (\Lambda s(a)\otimes t(a)\Lambda)_\Lambda \longrightarrow (\sfTr\Lambda)_\Lambda\longrightarrow 0.
\end{equation}
Moreover, as before, let $N$ (resp. $M$) be a left (resp. right) $\Lambda$-module.  Consider the $\Lambda$-bimodule  with external action $N\otimes M$. We have that $(N\otimes M)_\Lambda$ is isomorphic to $M\otimes_\Lambda N$ via the flip map.
Note that this is the degree $0$ instance of the graded isomorphism (see for example \cite[p.170 Corollary 4.4]{CARTANEILENBERG}):
$$\HH_*(\Lambda, N\otimes M)= \Tor_*^\Lambda(M,N).$$
  Thus,
$$(\Lambda x\otimes y\Lambda)_\Lambda= y\Lambda\otimes_\Lambda \Lambda x= y\Lambda x$$
which leads to the exact sequence
\begin{equation}\label{final coinvariants}
\oplus_{x\in Q_0}\ x\Lambda x \stackrel{f^*_\Lambda}{\longrightarrow} \oplus_{a\in Q_1} \ t(a)\Lambda s(a) \longrightarrow (\sfTr\Lambda)_\Lambda\longrightarrow 0.
\end{equation}
We underline that for $y,x\in Q_0$, the multiplicity of the vector space $y\Lambda x$ in the second  direct sum is the number of parallel arrows from $x$ to $y$.

Another easy and rather meticulous  computation gives a formula for $f^*_\Lambda$ in the context of (\ref{final coinvariants}). For $\lambda \in x\Lambda x$ we have
$$f^*_\Lambda(\lambda) = \sum_{ \substack{a\in Q_1\\t(a)=x}} \lambda a  -\ \sum_{ \substack{b\in Q_1\\s(b)=x}} b \lambda$$
where $\lambda a \in t(a)\Lambda s(a)$, that is the direct summand corresponding to $a$. Similarly, $b\lambda \in t(b)\Lambda s(b)$, that is the direct summand corresponding to $b$.

Let $C=\sum_{a\in Q_1}a \in \Lambda$. Note that for $\lambda\in \oplus_{x\in Q_0}\ x\Lambda x$ we have $$f^*_\Lambda(\lambda)= \lambda C-C\lambda.$$

To show that $\Ker f^*_\Lambda = Z\Lambda$, it is convenient as usual to consider the $k$-category $\C_\Lambda$ associated to $\Lambda$: its set of objects is $Q_0$, while the set of morphisms ${}_v\C_u$ from $u$ to $v$ is $v\Lambda u$; composition is given by the product of $\Lambda$. The center of $\Lambda$ viewed in this category is
$$\{({}_x\lambda_x)_{x\in Q_0} \ | \  {}_v\lambda_v  \  {}_v\alpha_u  =  {}_v\alpha_u\  \ {}_u\lambda_u\  \mbox{ for all } {}_v\alpha_u\in {}_v\C_u  \}.$$
On the other hand as already observed, in case of parallel arrows there is one direct summand for each arrow in $\oplus_{a\in Q_1} t(a)\Lambda s(a) $. Note also that $Q_0\cup Q_1$ is a set of generators of $\C_\Lambda$ as an algebra. Using these three observations, the proof of $\Ker f^*_\Lambda = Z\Lambda$ is immediate.\qed
\end{proof}

\section{\sf Hereditary, radical square zero and triangular monomial algebras}\label{cases}

In this section we compute the excess (see Definition \ref{excess}) of some families of algebras.

\subsection{\sf Hereditary algebras and algebras with the $\HH^2$ cancellation properties}

We first prove that the excess is zero for hereditary algebras. The  proof is based on the fact that the set of morphisms which do not factor through injectives is zero and we believe it provides a useful method in other contexts.

  Later in Theorem \ref{H^2s zero} we generalize the result for algebras with the $\HH^2$ cancellation properties (see the Introduction for the definition). Its proof relies  on the fact that for an algebra $\Lambda$ with the $\HH^2$ cancellation properties a formula for the dimension of $\HHH^1(\Lambda)$ is known, see \cite{CIBILS1988}.

\begin{theo}\label{hereditary}
Let $Q$ be a finite connected quiver without oriented cycles. Let $\Lambda=kQ$ be the corresponding hereditary algebra.
We have $e( \Lambda)=0.$
\end{theo}

\begin{proof}

We will show that  if $I$ is an injective $\Lambda$-bimodule, then $\Hom_{\Lambda-\Lambda}(I, \tau\Lambda)=0.$  \emph{A fortiori} $\mathcal{I}\Hom_{\Lambda-\Lambda} (\Lambda, \tau\Lambda)=0$. By Lemma \ref{excess is dim of factoring injectives}, it follows that $e(\Lambda)=0$ .

We have that $\mathsf{pd}_{\Lambda-   \Lambda} \Lambda \leq 1$. Indeed $kQ$ is the tensor algebra $T_{kQ_0} kQ_1$.
 It is well known (see for instance \cite[Theorem 2.3]{CIBILS1991})   that there is a  minimal  projective resolution of $kQ$ as a $kQ$-bimodule as follows:
\begin{equation}\label{projective resolution of kQ}
  0 \longrightarrow kQ\otimes_{kQ_0} kQ_1 \otimes_{kQ_0} kQ \longrightarrow  kQ\otimes_{kQ_0} kQ \longrightarrow kQ\longrightarrow 0.
  \end{equation}

 We recall \cite[Proposition 1.7(a) p. 319]{AUSLANDERREITENSMALO1995}: let $A$ be an algebra and let $M$ be  an indecomposable left $A$-module.  The projective dimension of $M$ is at most $1$   if and only if $\Hom_A (\sfD A, \tau M)=0$. We will use this result for $\Lambda$-bimodules, that is replacing $A$ by the enveloping algebra of $\Lambda$. We have supposed $Q$  connected, therefore $\Lambda$ is indecomposable as a $\Lambda$-bimodule, and the aforementioned proposition of \cite{AUSLANDERREITENSMALO1995} applies.

It follows that $\Hom_{\Lambda-\Lambda}(\sfD (\Lambda\otimes \Lambda), \tau \Lambda)=0$. Of course, for an algebra $A$, every  injective left $A$-module is isomorphic to a direct summand of a direct sum of copies of $\sfD A$, where $A$ is viewed as a right $A$-module and $A$ a left $A$-module. \qed
\end{proof}

\begin{coro}\cite{CIBILS1988,HAPPEL1988,CIBILS2000}
Let $B$ the set of paths of $Q$, and let $|yBx|$ be the number of paths from $x$ to $y$.
For $\Lambda= kQ$ we have
 $$\dim_k \HHH^1(\Lambda) = 1 - |Q_0| + \sum_{a\in Q_1}|t(a)Bs(a)| = \dim_k\tauh.$$
\end{coro}

We provide in the following a generalisation of Theorem \ref{hereditary} for algebras having the $\HH^2$ cancellation properties.

For a bound quiver algebra $\Lambda$ with the $\HH^2$ cancellation properties, the dimension of $\HHH^1(\Lambda)$ is known by \cite[p. 647]{CIBILS1988}. This allows to prove the following

\begin{theo} \label{H^2s zero}
The excess of a bound quiver algebra $\Lambda=kQ/I$ with the $\HH^2$ cancellation properties is zero.
\end{theo}

\begin{proof}
Let $B$ be the basis of paths of a bound quiver algebra.

We know from  \cite{CIBILS1988} that
\begin{equation*}
  \dim_k \HHH^1(\Lambda)=\dim_kZ\Lambda -\sum_{x\in Q_0} |xBx| + \sum_{x,y \in Q_0}|yBx||yQ_1x|.
\end{equation*}

Clearly $|yBx| = \dim_k y\Lambda x$. Hence by Theorem \ref{dimension of tau tilting cohomology} the equality of dimensions holds. \qed
\end{proof}
\begin{lemm}
A hereditary algebra $kQ$ has the $\HH^2$ cancellation properties.
\end{lemm}

\begin{proof}
It follows from (\ref{projective resolution of kQ}) that $\mathsf{pd}_{kQ-  kQ} kQ \leq 1$. Then  for any $kQ$-bimodule $X$ we have $\HH^2(kQ, X)= 0$.\qed
\end{proof}
\begin{rema}
We will show in Subsection \ref{subsection r^2=0} that not only the hereditary algebras have the $\HH^2$ cancellation properties.
\end{rema}

\subsection{\sf Radical square zero algebras}\label{subsection r^2=0}

A \emph{radical square zero algebra} is a bound quiver algebra of the form $kQ/F^2$.

Let $P$ and $P'$ be two sets of paths of a quiver $Q$. The set of \emph{parallel paths} is
$$P/\!/ P'=\{(p,p')\in P\times P' \ | \ s(p)=s(p') \mbox{ and } t(p)=t(p')\}.$$
For instance $Q_1/\!/ Q_0$  corresponds to the set of loops $Q_1^l=\{a\in Q_1 \ |\ s(a)=t(a)\}.$ 

We denote by $Q_i$ the set of paths of length $i$.

A $c$-crown is a quiver $C$ with $c$ vertices cyclically labelled and $c$ arrows, each one joining each vertex with the next one in the cyclic labelling. For instance a $1$-crown is a loop, and a $2$-crown is a two-way quiver $\cdot\leftrightarrows\cdot$.
The behaviour of the Hochschild cohomology of $kC/F^2$ is exceptional, see \cite{CIBILS1998} and it will be considered separately.

\begin{prop}
Let $Q$ be a connected quiver which is not a crown. The radical square zero algebra $\Lambda=kQ/F^2$ has the $\HH^2$ cancellation properties if and only if $Q_2/\!/ Q_1=\emptyset$.
\end{prop}

\begin{proof}
Since $r$ is a semisimple $\Lambda$-bimodule, the complex of cochains of Section 2 of \cite{CIBILS1998} has zero coboundaries and $\dim_kH^2(\Lambda, r)=|Q_2/\!/ Q_1|$.

Consequently if $\Lambda=kQ/F^2$ has the $\HH^2$ cancellation properties, then $|Q_2/\!/ Q_1|=0$.

Reciprocally, note first that if $|Q_2/\!/ Q_1|=0$ then $|Q_1/\!/ Q_0|=0$. From \cite[Theorem 3.1]{CIBILS1998} we have $$\dim_k \HHH^2(\Lambda) = |Q_2/\!/ Q_1| - |Q_1/\!/ Q_0| =0.$$
Hence if  $Q_2/\!/ Q_1=\emptyset$ then $\HH^2(\Lambda,r)=0=\HHH^2(\Lambda).$\qed
\end{proof}

There are zero excess algebras which does not have the $\HH^2$ cancellation properties,  see Corollary \ref{zero excess without H2 cancellation}. We first compute the excess of a radical square zero algebra.

\begin{prop}
Let $Q$ be a connected quiver which is not a crown and let $\Lambda=kQ/F^2$. We have $e(\Lambda)= |Q_1/\!/ Q_0|$.
\end{prop}

\begin{proof}
 We will use the formula of Theorem \ref{dimension of tau tilting cohomology}. First we observe the following:
\begin{itemize}
  \item $\dim_k Z\Lambda = 1 + |Q_1/\!/ Q_0|$.
  \item For $x\in Q_0$ we have $\dim_k x\Lambda x= 1 + |Q_1 /\!/ \{x\}|$. Therefore 
  $$\sum_{x\in Q_0}\dim_k x\Lambda x = |Q_0|+ |Q_1/\!/ Q_0|.$$
  \item 
  For $a\in Q_1 \backslash Q_1^l$ we have $\dim_k t(a)\Lambda s(a)= |t(a)Q_1s(a)|.$ 
  Then $$\sum_{a\in Q_1 \backslash Q_1^l} \dim_k t(a)\Lambda s(a) = | (Q_1 \backslash Q_1^l)  /\!/  (Q_1 \backslash Q_1^l)|.$$
  \item We have 
  $\sum_{a\in Q_1^l} \dim_k s(a)\Lambda s(a) = |Q_1^l| + |Q_1^l /\!/ Q_1^l|$ since the sum is over the loops -- not over the vertices which have loops.
  \item Finally
  $$\sum_{a\in Q_1} \dim_k t(a)\Lambda s(a)= |Q_1 /\!/ Q_0| + |Q_1 /\!/ Q_1|.$$
\end{itemize}

\begin{align*}
\displaystyle 
\dim_k\tauh =& 1 + |Q_1/\!/ Q_0| - |Q_0| - |Q_1/\!/ Q_0|  + |Q_1/\!/ Q_0|  + |Q_1/\!/ Q_1|  \\=& 
 1  - |Q_0|   + |Q_1/\!/ Q_0| + |Q_1/\!/ Q_1|.
\end{align*}
On the other hand we know from  \cite[Theorem 3.1]{CIBILS1998}, together with the observation  in the next paragraph,  that the following formula holds $$\dim_k \HHH^1(\Lambda) =  1  - |Q_0|  + |Q_1/\!/ Q_1|.$$\qed
\end{proof}

In the proof of Theorem 3.1 in \cite{CIBILS1998} it is stated that ``$D$ is injective for a positive $n$".  This is right for $n>0$.  However for $n=0$ the kernel of $D$ has dimension one. Hence the formula for $\dim_k \HHH^1(\Lambda)$ in the statement of \cite[Theorem 3.1]{CIBILS1998} has to be modified  by adding $1$.

\begin{coro}\label{zero excess without H2 cancellation}
Let $Q$ be a connected quiver which is not a crown. If the quiver has no loops and $Q_2/\!/ Q_1 \neq\emptyset$, then the radical square zero algebra $kQ/F^2$ does not have the $\HH^2$ cancellation properties and has zero excess.
\end{coro}

\begin{prop} Let $C$ be a $c$-crown, and let $\Lambda=kC/F^2$.

\begin{itemize}

\item If $c>1$, then $e(\Lambda) = 0,$

\item If $c=1$ and the characteristic of $k$ is not $2$, then $e(\Lambda)= 1,$

\item If $c=1$ and  the characteristic of $k$ is $2$, then $e(\Lambda)= 0.$
\end{itemize}

\end{prop}
\begin{proof} Observe that if $c>1$,  by  Theorem \ref{dimension of tau tilting cohomology} we have
$$\dim_k\HHH^1_\tau(kC/F^2) = 1 - c + c =1.$$
If $c=1$ then $kC/F^2= k[x]/(x^2)$, the algebra of dual numbers and
$$\dim_k \HHH^1_\tau (k[x]/x^2) = 2 - 2 + 2 = 2.$$

On the other hand, it is easy to compute that if $c>1$ then $\dim_k \HHH^1(\Lambda) =1$ regardless the characteristic of $k$.

If the characteristic of $k$ is not $2$, then $\dim_k \HHH^1(k[x]/(x^2))=1$ while if the characteristic of $k$ is $2$, then  $\dim_k \HHH^1(k[x]/(x^2))=2.$ \qed

\end{proof}

\subsection{\sf Triangular monomial algebras}

A \emph{monomial algebra} is a bound quiver algebra $\Lambda=kQ/I$ where $I$ is generated by a minimal set of paths denoted by $Z$. The algebra $\Lambda$ is \emph{triangular}  if it is a quotient of a finite dimensional hereditary algebra $kQ$, that is  if $Q$ has no oriented cycles. We set the following.
\begin{itemize}
  \item
$B$ is the set of paths of $Q$ which do not contain any path of $Z$.  Note that $Z=\emptyset$ if and only if $B$ is the set of all the paths of $Q$. Moreover,  $B$ gives a basis of $\Lambda$.
\item
$(Q_1/\!/B)_u$ is the set of  pairs $(a,\epsilon)\in Q_1/\!/B $ such that for every $\gamma\in Z$, replacing each occurrence of $a$ in $\gamma$  by $\epsilon$, gives a path which is $0$ in $\Lambda$.
Note that $\{ (a,a) \ | \ a\in Q_1\} \subset (Q_1/\!/B)_u.$

\item
$(Q_1/\!/B)_{nu}= (Q_1/\!/B) \setminus (Q_1/\!/B)_u,$ that is the set of  pairs $(a,\epsilon)\in Q_1/\!/B $ such that there exists $\gamma\in Z$ where $a$ occurs, verifying that at least one of the replacements of $a$ in $\gamma$ by $\epsilon$, gives a non zero path in $\Lambda$.
\end{itemize}

\begin{theo}
Let $\Lambda=kQ/\langle Z \rangle$ be a triangular monomial algebra. We have
$$e(\Lambda)= |(Q_1/\!/B)_{nu}|.$$
\end{theo}
\begin{proof}
When $Q$ has no oriented cycles, the formula for the dimension of $\HHH^1(\Lambda)$ given in  \cite{CIBILSSAORIN2001} is as follows:
\begin{equation}\label{Hoch1 for monomial without cycles}
\dim_k \HHH^1(\Lambda) = \dim_k Z\Lambda - |Q_0| + |(Q_1/\!/B)_u|.
\end{equation}

Note that since $Q$ has no oriented cycles, for all $x\in Q_0$ we have $x\Lambda x =k.$ Hence Theorem \ref{dimension of tau tilting cohomology} gives
\begin{equation}\label{monomial tau Hoch1}
\dim_k \tauh = \dim_k Z\Lambda - |Q_0| + |(Q_1/\!/B)|.
\end{equation}\qed
 \end{proof}

\begin{theo}
 Let $Q$ be a connected quiver without oriented cycles and  let $\Lambda=kQ/\langle Z \rangle$ be a triangular monomial algebra. The following are equivalent:
\begin{enumerate}[label=(\arabic*),itemsep=5pt]
  \item\label{Hoch1=0} $\HHH^1(\Lambda)=0.$
    \item \label{tree}$Q$ is a tree.
  \item \label{tauHoch1=0} $\tauh=0.$

\end{enumerate}
\end{theo}

\begin{rema} The equivalence between \ref{Hoch1=0} and \ref{tree} is proved without the triangular hypothesis in \cite[Theorem 2.2]{BARDZELL MARCOS 1998}.
\end{rema}

\begin{proof}

For \ref{Hoch1=0} implies \ref{tree}, the formula (\ref{Hoch1 for monomial without cycles}) gives
$$ 1 - |Q_0| + |(Q_1/\!/B)_u| =0.$$
 We have $\{ (a,a) \ | \ a\in Q_1\} \subset (Q_1/\!/B)_u$   hence $ 1 - |Q_0| + |Q_1| \leq 0.$ The Euler characteristic of the underlying graph of $Q$ is $\chi(Q)= |Q_0| - |Q_1|$, hence  $\chi(Q)\geq 1$. Any finite graph has the homotopy type of a graph with $1$ vertex and $n$ loops, which fundamental group is free on $n=1-\chi(Q)$ generators. We infer $n\leq 0$, hence $n=0$ and $Q$ is a tree.

Concerning \ref{tree} implies \ref{tauHoch1=0}, since $Q$ is a tree we have $\chi(Q)= |Q_0| - |(Q_1|=1$.  On the other hand $(Q_1/\!/B)=\{(a,a)\ \mid \ a\in Q_1\}$, and we have $|Q_1/\!/B|=|(Q_1|$. The formula \ref{monomial tau Hoch1} gives $\tauh=0.$

The implication \ref{tauHoch1=0} $\Rightarrow$ \ref{Hoch1=0} follows from Lemma \ref{excess is dim of factoring injectives}. \qed

\end{proof}

\begin{exam}
Let $Q$ and $R$ respectively denote the quivers
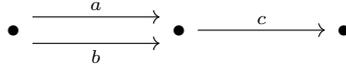
\[\begin{tikzcd}
	\bullet && \bullet && \bullet
	\arrow["b"', shift right=2, from=1-1, to=1-3]
	\arrow["a", shift left=2, from=1-1, to=1-3]
	\arrow["c", from=1-3, to=1-5]
\end{tikzcd}\]
and
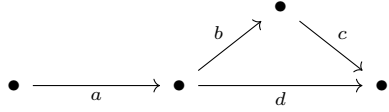
\[\begin{tikzcd}[ampersand replacement=\&]
	\&\&\& \bullet \\
	\bullet \&\& \bullet \&\& \bullet
	\arrow["a"', from=2-1, to=2-3]
	\arrow["b", from=2-3, to=1-4]
	\arrow["c", from=1-4, to=2-5]
	\arrow["d"', from=2-3, to=2-5]
\end{tikzcd}\]

 The following table lists the results for the corresponding monomial algebras:
 \begin{table}[ht]
\begin{tabular}{|c|c|c|c|c|c|}
\hline
\sf{Quiver} & Z        & $(Q_1/\!/B)_{nu}$ & $e(\Lambda)$ & $\dim_k \HHH^1(\Lambda)$ & $\dim_k\tauh$ \\ \hline
$Q$    & $\{ca\}$ & $\{(a,b)\}$       & $1$      & $2$                & $3$           \\ \hline
$R$    & $\{ba\}$ & $\emptyset$       & $0$      & $2$                & $2$           \\ \hline
$R$    & $\{da\}$ & $\{(d,bc)\}$      & $1$      & $1$                & $2$           \\ \hline
\end{tabular}
\end{table}

\end{exam}

\section{\sf The excess}\label{The excess}

The proof of the next result relies on the calculation of the dimensions of two vector spaces and the observation that they are equal. An explicit isomorphism between these two vector spaces remains unknown to us.

\begin{prop}
Let $\Lambda=kQ/I$ be a bound quiver algebra. We have
$$\dim_k \HH^1(kQ, kQ/I) = \dim_k \tauh.$$
\end{prop}
\begin{proof}
Let $X$ be a $kQ$-bimodule. We assert that
\begin{equation}\label{H1(kQ,X)}
\dim_k \HH^1(kQ,X)= \dim_k X^{kQ} - \sum_{x\in Q_0} \dim_k xXx + \sum_{a\in Q_1} \dim_k t(a)Xs(a).
\end{equation}
Recall the projective resolution of $kQ$ as a $kQ$-bimodule (\ref{projective resolution of kQ})
\begin{equation}
  0 \longrightarrow kQ\otimes_{kQ_0} kQ_1 \otimes_{kQ_0} kQ \stackrel{g}{\longrightarrow}  kQ\otimes_{kQ_0} kQ \longrightarrow kQ\longrightarrow 0.
  \end{equation}
The functor $\Hom_{kQ-kQ}(-,X)$ gives the complex of cochains
\begin{equation*}
  0\longrightarrow \Hom_{kQ-kQ}(kQ\otimes_E kQ, X) \stackrel{g^*}{\longrightarrow} \Hom_{kQ-kQ}(kQ\otimes_E kQ_1\otimes_E kQ, X)\longrightarrow 0
\end{equation*}
where $\Ker g^* = \HH^0(kQ, X)$ and $\Coker g^* = \HH^1(kQ, X)$.
The same way we have obtained (\ref{exact Tr E-E}) and (\ref{quasi final Tr}) leads to an exact sequence
\begin{equation*}
  0\longrightarrow \HH^0(kQ, X)\longrightarrow \oplus_{x\in Q_0} \ xXx \stackrel{g^*}{\longrightarrow} \oplus_{a\in Q_1} \ t(a)Xs(a) \longrightarrow \HH^1(kQ, X) \longrightarrow 0.
\end{equation*}
which gives the equality (\ref{H1(kQ,X)}).

We assert that if $X$ is a $kQ/I$-bimodule, that is $IX=XI=0$, then $X^{kQ}= X^{kQ/I}$. Indeed, we have
$$X^{kQ}= \Hom_{kQ-kQ}(kQ, X)=  \Hom_{kQ/I-kQ/I}(kQ/I, X)= X^{kQ/I}.$$
Note that for $X=\Lambda$ we have $\Lambda^\Lambda= Z\Lambda$.  We obtain the following
$$\dim_k \HH^1(kQ,kQ/I)= \dim_k Z\Lambda - \sum_{x\in Q_0} \dim_k x\Lambda x + \sum_{a\in Q_1} \dim_k t(a)\Lambda s(a)$$
which is the same formula than the one for $\dim_k\tauh$ in Theorem \ref{dimension of tau tilting cohomology}.\qed
\end{proof}

Next we recall Corollary 2.4 of \cite{CIBILS1991}.

\begin{prop}\cite{CIBILS1991}\label{4}
Let $\Lambda=kQ/I$ be a bound quiver algebra and let $X$ be a $\Lambda$-bimodule. There is an exact sequence
\begin{equation*}
0\longrightarrow \HH^1(\Lambda, X)\longrightarrow \HH^1(kQ,X)\longrightarrow \Hom_{kQ-kQ}(I/I^2, X)\longrightarrow \HH^2(\Lambda, X)\longrightarrow 0.
\end{equation*}
\end{prop}

An immediate consequence of the above is a formula for the excess of an algebra, which involves the dimension of the Hochschild cohomology in degree $2$.
\begin{theo}\label{exceso}
Let $\Lambda=kQ/I$ be a bound quiver algebra. We have
$$e(\Lambda) = \dim_k \Hom_{kQ-kQ}(I/I^2, \Lambda)\ -\ \dim_k \HHH^2(\Lambda).$$
\end{theo}

\begin{rema}
If $I=0$, then $\HHH^2(kQ)=0$ and so $e(kQ)=0$. This confirms Theorem \ref{hereditary}, as well as Theorem \ref{H^2s zero} for an hereditary algebra.
\end{rema}

We infer  three    corollaries for a bound quiver algebra $\Lambda=kQ/I$.

\begin{coro}
If $\Lambda$ verifies the $\HH^2$ cancelling properties then $$\Hom_{kQ-kQ}(I/I^2, \Lambda)=0.$$
\end{coro}

\begin{coro} \label{H2 is zero}
If $e(\Lambda)=0$ then
$$ \HHH^2(\Lambda) = \Hom_{kQ-kQ}(I/I^2, \Lambda).$$
\end{coro}

\begin{coro}
The algebra $\Lambda$ is $\tau$-rigid as a $\Lambda$-bimodule if and only if
\begin{itemize}
  \item $\HHH^1(\Lambda)=0,$
  \item $\dim_k \HHH^2(\Lambda)= \dim_k \Hom_{kQ-kQ}(I/I^2, \Lambda).$
\end{itemize}
\end{coro}

\footnotesize
\noindent C.C.:\\
Institut Montpelli\'{e}rain Alexander Grothendieck, CNRS, Univ. Montpellier, France.\\
{\tt Claude.Cibils@umontpellier.fr}

\medskip
\noindent M.L.:\\
Instituto de Matem\'atica y Estad\'\i stica  ``Rafael Laguardia'', Facultad de Ingenier\'\i a, Universidad de la Rep\'ublica, Uruguay.\\
{\tt marclan@fing.edu.uy}

\medskip
\noindent E.N.M.:\\
Departamento de Matem\'atica, IME-USP, Universidade de S\~ao Paulo, Brazil.\\
{\tt enmarcos@ime.usp.br}

\medskip
\noindent A.S.:
\\IMAS-CONICET y Departamento de Matem\'atica,
Facultad de Ciencias Exactas y Naturales,\\
Universidad de Buenos Aires, Argentina. \\{\tt asolotar@dm.uba.ar}

\end{document}